# One-vs-one Threat-Aware Weaponeering with Basic Engagement Zones

Alexander Von Moll[1], Dejan Milutinović[2], Isaac Weintraub[1], and David W. Casbeer[1]

*Abstract*—In this paper we address the problem of 'weaponeering', i.e., placing the weapon engagement zone (WEZ) of a vehicle on a moving target, while simultaneously avoiding the target's WEZ. A WEZ describes the lethality region of a range-limited weapon considering both the range of the weapon along with the state of the target. The weapons are assumed to have simple motion, while the vehicles carrying the weapons are modeled with Dubins dynamics. Three scenarios are investigated and are differentiated in the assumptions that can be made about the target in the process of the vehicle control design: 1) no knowledge of target control, 2) avoid unsafe positions assuming the target's optimal control, 3) full knowledge of target's optimal control. The engagement is formulated as a stochastic optimal control problem with uncertainty in the target's control modeled using a noise parameter applied to the target's control input. After discretizing the Hamilton-Jacobi-Bellman equation, Value iteration is then used to obtain an approximate solution for the optimal vehicle control and time-to-go. Simulation results support usage of the first paradigm: assume no knowledge of the target's control.

## I. INTRODUCTION

Pursuit of a moving target is a common component of real-world missions ranging from search and rescue [1], rendezvous [2], and missile intercept [3]. In fact, pursuit-evasion scenarios were among the first investigated as part of the genesis of the theory of differential games [4], [5]. A list of many classical pursuit-evasion differential game scenarios is given in [6].

Much of the literature on pursuit-evasion is concerned with point capture wherein one agent is said to capture the other if they are spatially collocated. Aside from the mathematical issues with this model, a more realistic capture criterion is the finite capture radius (c.f., e.g., [7]). There, capture is said to occur once the distance between the agents is within the specified capture radius. Physically speaking, this could correspond to a blast radius, where the blast propagates very fast. Depending on the application (e.g., vehicles as opposed to missiles), an even more realistic model is one in which the pursuer must get within some reasonable range of its target whereupon it must fire a weapon to actually intercept the target. Furthermore, such a model allows for the consideration of the evader/target's defensive capabilities by assuming it has its own weapon onboard. The 'reasonable range' that must be reached in order for the weapon to have a chance at reaching its target corresponds to the weapon engagement zone (WEZ).

However, the WEZ-capture paradigm is in stark contrast to the point capture and radius capture models wherein the evader is essentially helpless. There have been several recent works concerning WEZs. One possible definition of a WEZ is given by the following [8]:

"Given a Mobile Agent, $A$, Threat, $T$, their respective dynamic models, $\dot{\boldsymbol{x}}_A$, $\dot{\boldsymbol{x}}_T$, and an Agent strategy, $u_A(t)$, a [Weapon] Engagement Zone (WEZ) is a region of the state space in which it is possible for the Threat to neutralize the Mobile Agent if the latter does not deviate from its current strategy."

Most WEZ-related works have had to do with planning a path for a vehicle to reach a goal while avoiding entering the WEZ; indeed that application was revisited in [8], though its main purpose was to define and derive WEZs based on basic, first principles models. The important characteristic of WEZs is that they depend on the relative position of the vehicles and their headings, which distinguishes these works from the large body of research on path planning in the presence of obstacles. In the case of WEZs based on kinetic weapons, the WEZ accounts for the weapon time of flight as well as the vehicle's course. In essence, the weapon time of flight defines a shift while the vehicle's course defines the direction of the shift. Prior works considered a notional WEZ region without an underlying threat model.

### A. Literature Review

Weintraub et al. also considered optimal navigation of a single vehicle around a single threat to a specified goal region [9]. It is important to note that the WEZ location (i.e., the origin of the dynamic region, about which it rotates and scales) was assumed to be static (thus representing, for example, a stationary, ground-based defensive site). The solution was obtained using pseudospectral collocation. Using a similar approach, the solutions were extended to a single vehicle navigating around and between two WEZs [10]. Some challenges were uncovered in the process of extending to two WEZs: namely, the need for a high quality initial guess, the fact that there are several distinctly different 'options' (going around versus between), and scalability of the solution approach.

Recently, several works have extended the scenario to consider navigation of a single agent among many (i.e., > 2) WEZs to a specified goal. In order to address the scalability challenges observed previously, different solution approaches were taken. Wolek et al. utilized an rapidly-exploring random tree (RRT) approach which are able to find feasible path plans quickly and continuously refine the plan as the algorithm runs [11]; this approach was able to plan paths around 24 WEZs without issue. The main advantage of the RRT-based approach is near-real time execution (i.e, path plans on the order of several seconds). However, the only guarantee with regard to optimality is that the RRT* algorithm used converges to the true optimal path in the limit as time goes to infinity. Milutinović et al., on the other hand, utilized stochastic optimal control, discretization, and Value iteration to obtain an optimal feedback controller that could handle noise (e.g., due to wind) [12]. This method was far more computationally intensive (suitable for pre-mission planning only), and yielded a policy which provides the optimal control for any point in the state space. Three WEZs were considered in [12], however the method is suitable for many more WEZs.

Most of these works have focused on optimal navigation of a vehicle around a stationary threat. However, a new thread has emerged, however, concerning the so-called weaponeering problem which entails a vehicle with a WEZ of its own seeking to place its WEZ on a moving target in minimum time. Perhaps the simplest case is when the target moves with constant speed and course, the vehicle moves with simple motion (single integrator kinematics), and the WEZ is modeled by the range-limited pursuit-evasion basic engagement zone (BEZ) from [8]. This case is treated by Weintraub et al. in [13]. These works are preliminary in the sense that simplified motion models are used and it is assumed that the target does not maneuver. Moreover, there is no consideration for avoidance of the target's WEZ (if it were to have one).

### B. Proposed Approach

The goal is to address a "WEZ versus WEZ" scenario in which the vehicle 'wins', in a sense, if it can get the target in its WEZ and 'loses' if it gets into the target's WEZ. As an improvement over the previous works on weaponeering we will consider both vehicles to be Dubins vehicles (i.e., constant speed, turn-rate-constrained vehicles representative of fixed wing aircraft).

Conceptually, the scenario is nearly identical to the scenario explored in [14] wherein the vehicle desired to get in the tail of its target with heading roughly aligned with the target heading.

From a differential game-theoretic standpoint, similar scenarios have been considered, often under the title *Game of Two Cars*, e.g., [15], [16]. Merz also devoted specific attention to the all important question of whether or not the vehicle should be pursuing or evading based on the current configuration [17].

Like [13], we will at first consider the WEZs to be modeled by a BEZ from [8]. However, the proposed methodology will be extensible to other models. The proposed solution methodology is based on stochastic optimal control and will follow a similar process as presented in [14].

### C. Contributions and Outline

The contributions of this paper are that it is the first work to:
1) incorporate avoidance of a moving WEZ
2) consider placement of a WEZ onto a maneuvering target (since [13] considered a non-maneuvering target)
3) solve for an optimal feedback controller for placement of a WEZ onto a target

The controllers obtained in this work are suitable for real-time implementation onboard fixed-wing aircraft flying at fixed altitudes and constant speeds. Additionally the solution framework could be applied to a variety of scenarios – the target's WEZ could be replaced with a simple keep-out range, for example.

The remainder of the paper is outlined as follows. Section II provides the details of the mathematical formulation. Section IV describes some practical matters pertaining to the computation of the solution (which is relatively burdensome). Section V details the assumptions embedded in the three aforementioned scenarios: the baseline (no target control knowledge), avoidance, and adversarial (assume target uses its optimal control). Section VI contains numerical results and the paper is concluded in Section VII.

## II. PROBLEM FORMULATION

Define two vehicles, vehicle $A$ (for "Agent") and vehicle $T$ (for "Target"), both with Dubins vehicle dynamics. We seek to solve for $A$'s control in this scenario and we assume that $T$'s control is unknown. Therefore, $A$'s dynamics are given by

$$\begin{aligned} \mathrm{d}x_A &= v_A \cos(\theta_A)\,\mathrm{d}t \\ \mathrm{d}y_A &= v_A \sin(\theta_A)\,\mathrm{d}t \\ \mathrm{d}\theta_A &= u_A\,\mathrm{d}t \end{aligned} \quad (1)$$

where $v_A$ is the Agent's constant speed and $u_A \in [-\bar{u}_A, \bar{u}_A]$ is its control. Meanwhile, $T$'s control $u_T \in [-\bar{u}_T, \bar{u}_T]$ is assumed to be unknown to $A$. Therefore, $T$'s heading increment is approximated using some nominal control $u_T$ along with some noise:

$$\begin{aligned} \mathrm{d}x_T &= v_T \cos(\theta_T)\,\mathrm{d}t \\ \mathrm{d}y_T &= v_T \sin(\theta_T)\,\mathrm{d}t \\ \mathrm{d}\theta_T &= u_T\,\mathrm{d}t + \sigma\,\mathrm{d}w \end{aligned} \quad (2)$$

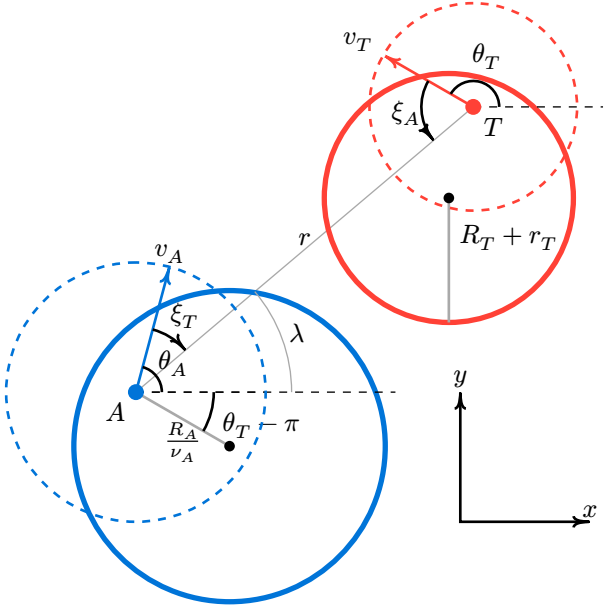

Fig. 1. Schematic depiction of the WEZ versus WEZ scenario.

where $v_T$ is the Target's constant speed, $\sigma$ is a noise intensity parameter, and $dw$ are unit intensity Wiener process increments.

The above system of equations is comprised of 6 states; therefore, it is advantageous to derive a reduced set of states. The aspect angles $\xi_A$ and $\xi_T$ are defined as follows and are shown in Fig. 1,

$$\xi_A = \lambda - \theta_T + \pi, \quad \xi_T = \lambda - \theta_A \quad (3)$$

where the bearing angle of the Target w.r.t. the Agent is

$$\lambda = \operatorname{atan}(y_T - y_A, x_T - x_A). \quad (4)$$

For example, $\xi_A$ can be thought of as "the aspect angle of $T$ as seen by $A$". Furthermore, let $r$ be the distance from $A$ to $T$,

$$r = \sqrt{(x_T - x_A)^2 + (y_T - y_A)^2}. \quad (5)$$

The increment in $r$ is obtained by taking the derivative of (5) w.r.t. the original states and simplifying,

$$dr = -(v_T \cos \xi_A + v_A \cos \xi_T) \, dt = b_r \, dt \quad (6)$$

Similarly, the increments in $\xi_A$ and $\xi_T$ are obtained by differentiating (3):

$$\begin{aligned}
d\xi_A &= \left[-u_T + \frac{1}{r}(v_T \sin \xi_A + v_A \sin \xi_T)\right] dt + \sigma \, dw \\
&= b_{\xi_A} \, dt + \sigma \, dw \\
d\xi_T &= \left[-u_A + \frac{1}{r}(v_T \sin \xi_A + v_A \sin \xi_T)\right] dt \\
&= b_{\xi_T} \, dt.
\end{aligned} \quad (7)$$

Let the state of the system be denoted $\mathbf{x} := (r, \xi_A, \xi_T)^\top$.

Both vehicles are assumed to carry onboard weapons modeled by points moving with simple motion, constant speed, with limited range, and a capture radius. This weapon model is the pursuit-evasion basic engagement zone described in [8]. For example, the Agent carries weapons with speed $\nu_A v_T$, where $\nu_A > 1$, range $R_A > 0$, and capture radius $r_A > 0$. The weapons carried by the Target are defined similarly.

The engagement ends when one of the vehicles enters the other vehicle's WEZ. For example, the set of states for which the Target is in the Agent's WEZ is given by

$$\mathcal{W}_A := \{\mathbf{x} \mid r \leq \rho_A(\xi_A)\}, \quad (8)$$

where $\overline{AT}$ is the distance between the vehicles. The set of states for which the Agent is in the Target's WEZ is similarly defined as

$$\mathcal{W}_T := \{\mathbf{x} \mid r \leq \rho_T(\xi_T)\}. \quad (9)$$

The quantities $\rho_A$ and $\rho_T$ represent, respectively, the distance from the vehicle to the boundary of its WEZ as a function of aspect angle:

$$\rho_\bullet(\xi_\bullet) = \frac{R_\bullet}{\nu_\bullet}\left[\cos \xi_\bullet + \sqrt{\cos^2 \xi_\bullet - 1 + \frac{(R_\bullet + r_\bullet)^2}{R_\bullet^2}}\right], \quad (10)$$

where $\bullet \in \{A, T\}$.

The goal for $A$ is to get its WEZ onto $T$ in minimum time while avoiding the WEZ of $T$. Thus the cost of $A$, which it wishes to minimize, is expressed as

$$J_A(u_A(\cdot)) = g(\mathbf{x}(t_f)) + \int_0^{t_f} dt, \quad (11)$$

where $t_f$ is the first time instant for which $T$ is in the WEZ of $A$, i.e., $\mathbf{x} \in \mathcal{W}_A$, or $A$ is in the WEZ of $T$, and

$$g(\mathbf{x}) = \begin{cases} 0, & \text{if } \mathbf{x} \in \mathcal{W}_A \\ M, & \text{if } \mathbf{x} \in \mathcal{W}_T. \end{cases} \quad (12)$$

Here, $M$ represents some large (positive) penalty for $A$ ending inside the WEZ of $T$. The goal is to obtain the optimal control $u_A^*$ which minimizes the cost, (11), in expectation.

### III. SOLUTION APPROACH

The overall solution approach utilized here is based upon the numerical stochastic optimal control techniques compiled in [18] and demonstrated in [14], [19], [20]. In particular, [14] is most similar to the present work as that work considered a 1 vs. 1 engagement of two vehicles modeled with Dubins car kinematics.

In general, the optimal control, $u_A^*$, may be obtained by solving the Hamilton-Jacobi-Bellman (HJB) partial differential equation [18], [21]:

$$0 = \min_{u_A}\{\mathcal{L}V(\mathbf{x}) + 1\}, \quad (13)$$

where $V$ is the Value function which describes the optimal cost of going from the state $\mathbf{x}$ to termination, and $\mathcal{L}$ is the differential operator

$$\mathcal{L} = b_r \frac{\partial}{\partial r} + b_{\xi_A} \frac{\partial}{\partial \xi_A} + b_{\xi_T} \frac{\partial}{\partial \xi_T} + \frac{1}{2}\sigma^2 \frac{\partial^2}{\partial \xi_A^2} \quad (14)$$

In practice it can be quite difficult to obtain an analytic solution to the HJB, except for very simple problems. Therefore, the methodology from [18] is based upon discretizing the HJB equation using a locally consistent Markov chain approximation to compute the optimal control. The Markov chain approximation is locally consistent in the sense that the mean and variance converge to the mean and variance of the continuous time process. Within this Markov chain, the transition probabilities are dependent on the control of both vehicles. The discretization essentially allows one to solve an approximate version of (13) via dynamic programming using Value iterations. As a result, one obtains $V^h$ and $u_A^{*h}$ which are discretized approximations of the Value function and associated optimal control with respect to the discretization $h$.

First, the HJB is discretized in the state space. Let $\Delta r$, $\Delta \xi_A$, and $\Delta \xi_T$ be the step sizes in each of the dimensions. Let $\mathbf{x}^h := (r^h, \xi_A^h, \xi_T^h)^\top$ be a particular state within the discretization and define $V^h := V(\mathbf{x}^h)$. Also, define $\underline{r^h} := r^h - \Delta r$ and $\overline{r^h} := r^h + \Delta r$ with $\underline{\xi_A}$, $\overline{\xi_A}$, $\underline{\xi_T}$, and $\overline{\xi_T}$ defined similarly. Then (13) is approximated with the following (i.e., the upwind derivative approximations),

$$b_r \frac{\partial V}{\partial r} \approx \frac{b_{r^h}^+}{\Delta r}\left(V(\overline{r^h}, \xi_A^h, \xi_T^h) - V^h\right)$$
$$- \frac{b_{r^h}^-}{\Delta r}\left(V^h - V(\underline{r^h}, \xi_A^h, \xi_T^h)\right)$$
$$b_{\xi_A} \frac{\partial V}{\partial \xi_A} \approx \frac{b_{\xi_A^h}^+}{\Delta \xi_A}\left(V(r^h, \overline{\xi_A^h}, \xi_T^h) - V^h\right)$$
$$- \frac{b_{\xi_A^h}^-}{\Delta \xi_A}\left(V^h - V(r^h, \underline{\xi_A^h}, \xi_T^h)\right)$$
$$b_{\xi_T} \frac{\partial V}{\partial \xi_T} \approx \frac{b_{\xi_T^h}^+}{\Delta \xi_T}\left(V(r^h, \xi_A^h, \overline{\xi_T^h}) - V^h\right) \quad (15)$$
$$- \frac{b_{\xi_T^h}^-}{\Delta \xi_T}\left(V^h - V(r^h, \xi_A^h, \underline{\xi_T^h})\right)$$
$$\frac{1}{2}\sigma^2 \frac{\partial^2 V}{\partial \xi_A^2} \approx \frac{\sigma^2}{2(\Delta \xi_A)^2}\left(V(r^h, \overline{\xi_A^h}, \xi_T^h) - V^h\right)$$
$$- \frac{\sigma^2}{2(\Delta \xi_A)^2}\left(V^h - V(r^h, \underline{\xi_A^h}, \xi_T^h)\right),$$

where $b_{r^h}^+ = \max\{0, b_{r^h}\}$, $b_{r^h}^- = \max\{0, -b_{r^h}\}$ and $b_{\xi_A^h}^+$, $b_{\xi_A^h}^-$, $b_{\xi_T^h}^+$, and $b_{\xi_T^h}^-$ are defined similarly. Since the state space has been discretized, the time step depends on the other parameters, i.e., let $\Delta t$ be the time implicit time step [18, pp. 112]:

$$\Delta t = \left(\frac{|b_{r^h}|}{\Delta r} + \frac{|b_{\xi_A^h}|}{\Delta \xi_A} + \frac{|b_{\xi_T^h}|}{\Delta \xi_T} + \frac{\sigma^2}{(\Delta \xi_A)^2}\right)^{-1} \quad (16)$$

where $|b_\bullet| := b_\bullet^+ + b_\bullet^-$ for $\bullet \in \{r, \xi_A, \xi_T\}$. Substituting the upwind derivative approximations, (15), and the time step, (16), into (13) and moving all the terms that include $V^h$ to the left side of the expression yields

$$V^h = \min_{u_A}\Big\{\Delta t + p_{\Delta r^+} V\left(\overline{r^h}, \xi_A^h, \xi_T^h\right) +$$
$$+ p_{\Delta r^-} V(\underline{r^h}, \xi_A^h, \xi_T^h) + p_{\Delta \xi_A^+} V\left(r^h, \overline{\xi_A^h}, \xi_T^h\right)$$
$$+ p_{\Delta \xi_A^-} V\left(r^h, \underline{\xi_A^h}, \xi_T^h\right) + p_{\Delta \xi_T^+} V\left(r^h, \xi_A^h, \overline{\xi_T^h}\right) \quad (17)$$
$$+ p_{\Delta \xi_T^-} V\left(r^h, \xi_A^h, \underline{\xi_T^h}\right)\Big\},$$

where the quantities

$$p_{\Delta r^\pm} = \Delta t\, b_{r^h}^\pm / \Delta r$$
$$p_{\Delta \xi_A^\pm} = \Delta t\left(b_{\xi_A^h}^\pm / \Delta \xi_A + \sigma^2/(2(\Delta \xi_A)^2)\right) \quad (18)$$
$$p_{\Delta \xi_T^\pm} = \Delta t\, b_{\xi_T^h}^\pm / \Delta \xi_T$$

represent the transition probabilities of the Markov chain from cells neighboring $\mathbf{x}^h$.

Value iterations [22] may be performed on $V^h$ using (17) starting from an initial guess for all of the cells in the discretization (see Algorithm 1). The computational domain is defined as

$$\Omega = [0, r_{\max}] \times [-\pi, \pi - \xi_A] \times [-\pi, \pi - \Delta \xi_T]. \quad (19)$$

Because of the periodicity of the $\xi_A$ and $\xi_T$ states, the cells at $\xi_A = -\pi$ are neighbors with the cells at $\xi_A = \pi - \Delta \xi_A$, for example. Far from the Target, the policy should be essentially the same (i.e., the Agent should turn and point towards the Target, thereby entering the computational domain) Therefore, at $r = r_{\max}$, the reflective boundary condition is used, i.e., $V(r_{\max}, \xi_A^h, \xi_T^h) = V(r_{\max} - \Delta r, \xi_A^h, \xi_T^h)$. This helps to ensure that the stochastic process remains inside the computational domain [14]. In practice, one must choose $r_{\max}$ appropriately for the particular mission and, perhaps, specify some other behavior or control for states where $r > r_{\max}$. When the state is such that $A$ is inside the WEZ of $T$ the Value is $M$ (as shown in (12)). Similarly, when $\mathbf{x}^h \in \mathcal{W}_A$ the Value is 0.

---

**Algorithm 1:** Value Iteration

---

1 initialize $V(\mathbf{x}^h) = u_A^{*h} = 0$ for all $\mathbf{x}^h \in \Omega$
2 **for** $z \in 1, ..., N_{\text{iter}}$ **do**
3      **for** $\mathbf{x}^h \in \Omega$ **do**
4          **if** $r^h \in \{0, r_{\max}\}$ **then**
5              $V^h = V^h_{\text{neighbor}}$
6          **else if** $\mathbf{x}^h \in \mathcal{W}_T$ **then**
7              $V^h = M$
8          **else if** $\mathbf{x}^h \in \mathcal{W}_A$ **then**
9              $V^h = 0$
10          **else**
11              $V^h = \min_{u_A \in \{-\bar{u}_A, 0, \bar{u}_A\}}\{...\}$    eq. (17)
12              $u_A^{*h} = \arg\min_{u_A \in \{-\bar{u}_A, 0, \bar{u}_A\}}\{...\}$    eq. (17)

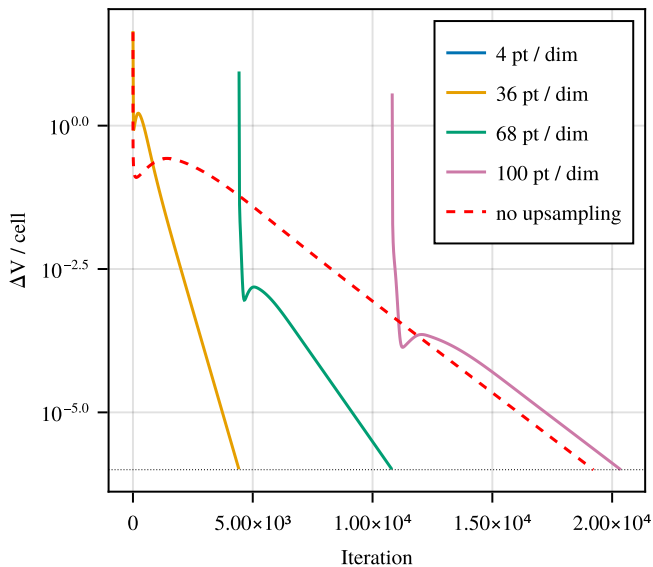

Fig. 2. Comparison of convergence of the Value iterations with and without upsampling. A speedup of 1.63x was achieved by utilizing upsampling.

## IV. COMPUTATIONAL CONSIDERATIONS

The fixed point of (17) is the solution to the discretized HJB. It is necessary to run Algorithm 1 for a large number of iterations, $N_{\text{iter}}$, to reach a reasonable level of convergence (measured in terms of the average change in Value) over the whole domain. Additionally, the state is 3-dimensional which entails significant computation; for example, performing 1,000 updates of a grid with 100 points in each dimension entails computing (17) on the order of one billion times. Fortunately, Algorithm 1 may be parallelized. For example, in Line 3 of Algorithm 1, the grid can be divided up among multiple threads. From a thread safety standpoint, one may use a copy of $V_{z-1}$ to evaluate the RHS of (17) so that no two threads attempt to read and write the same grid cell simultaneously.

Convergence can also be achieved faster by replacing Line 1 of Algorithm 1 with a more informative initialization. Upsampling is one useful technique for obtaining an informative initialization for the Value. First, Algorithm 1 is run for a relatively small grid size, the resulting grid of Values is then linearly interpolated over a finer grid, and Algorithm 1 is run again starting from these interpolated values. Fig. 2 shows the convergence history for sequential upsampling and hot-starting. Upsampling occurs at each of the discontinuities. The process reaches a desired threshold of convergence of $\Delta V = 1 \times 10^{-6}$ per cell within about the same total number of iterations as the 'no upsampling case'. This results in faster execution because the iterations performed on smaller grid sizes are significantly cheaper.

Note that the control input was discretized such that $u_A^{*h} \in \{-\bar{u}_A, 0, \bar{u}_A\}$ as shown in Algorithm 1, Line 12. For vehicles with Dubins dynamics, as in (1), such a discretization is reasonable because time-optimal paths are generally comprised of segments of maximum turn rate and straight lines (c.f. [23]).

Finally, the control policy resulting from Algorithm 1 is encoded as a lookup table, which may not be the most suitable for simulation, depending on the timestep. Therefore, linear interpolation is used when simulating the control policy $u_A^*$.

## V. ALTERNATIVE CONTROLLERS

In this section, we introduce *Baseline*, *Avoid*, and *Adversarial* controllers. The purpose in defining these controllers is to investigate whether there is a significant performance gain for $A$ against $T$ under different assumptions under which $A$'s controller is designed. Each of them are computed using Algorithm 1, though some require slight modification. Recall, in particular, that the derivations in the preceding sections were based on knowledge of a nominal controller $u_T$.

### A. Baseline Controller

Let $u_A^B(\mathbf{x})$ denote $A$'s Baseline controller. This controller is a numerical stochastic optimal controller for $A$ which assumes $u_T(\mathbf{x}) = 0$ for all $\mathbf{x} \in \Omega$ and $\sigma > 0$. This assumption reflects that $A$ does not know the $T$'s intention and $\sigma > 0$ completely encompasses uncertainty about changes in $T$'s heading. A small/large $\sigma$ value describes smaller/larger uncertainty, i.e., expectation that $T$ follows trajectories with smaller/larger changes in its heading. The Baseline controller Value function is $V_A^B$. In a similar way we can define $u_T^B$, a Baseline controller for $T$, which is computed by setting $\sigma = 0$ in (2) and replacing $u_A \, dt$ with $\sigma \, dw$ in (1). The corresponding Value function is $V_T^B$.

### B. Avoid Controller

The Avoid controller is based off of a concept from [14] wherein the blue vehicle (Agent) seeks to avoid states in which the red vehicle's (Target) time-to-go is smaller. The set of states for which $T$'s time-to-go is shorter is given by

$$\mathcal{D} = \{\mathbf{x} \in \Omega \mid V_T^B(\mathbf{x}) < V_A^B(\mathbf{x})\}. \tag{20}$$

Finally, $A$'s Avoid controller is obtained by assuming $u_T = 0$ and running Algorithm 1 with Line 6 modified as follows:

---
**Algorithm 2:** Value Iteration (mod. for Avoid controller)

---
6          **else if** $\mathbf{x}^h \in \mathcal{W}_T \cup \mathcal{D}$ **then**

---

In other words, the Value is set to $M$ when the state lies in the avoidance set $\mathcal{D}$.

### C. Adversarial Controller

Finally, the Adversarial controller is a controller for $A$ which assumes that $T$ *uses* its Baseline controller. This is a more aggressive assumption compared to the Avoid controller, which simply seeks to stay away from "dangerous" states in which $T$ could *potentially* capture $A$ first. The Adversarial controller assumes that the Target is actively trying to get $A$ into its WEZ (for all states in the domain). Note that neither the Baseline controller nor the Adversarial controller are necessarily the equilibrium policy associated with a two-player, stochastic zero-sum differential game as both controllers make particular assumptions about the Target's control which also may not correspond to $T$'s equilibrium policy.

TABLE I. VEHICLE PARAMETER SETTINGS FOR NUMERICAL RESULTS

|  | Agent | Target |
|---|---|---|
| speed | $v_A = 1$ | $v_T = 0.8$ |
| weapon speed ratio | $\nu_A = 1.2$ | $\nu_T = 1.1$ |
| weapon range | $R_A = 1$ | $R_T = 0.9$ |
| capture radius | $r_A = 0.2$ | $r_T = 0.15$ |
| max turn rate | $\bar{u}_A = 1$ | $\bar{u}_T = 1$ |

TABLE II. COMPUTATIONAL PARAMETERS

| number of iterations | $2 \times 10^4$ |
|---|---|
| convergence tolerance, $\Delta V$/cell | $1 \times 10^{-6}$ |
| points per dimension | 100 |
| max distance, $r_{max}$ | 10 |
| noise intensity, $\sigma$ | 1 |

## VI. RESULTS

In this section, the numerical results of the solution are presented along with results of simulations of the Agent kinematics, (1), under various controllers, and noise-free Target kinematics, (2), under various controllers. All of the controllers were computed with the vehicle parameters from Table I. $A$ has the advantage in all parameters except in turn rate. Unless otherwise specified, the controllers were all computed using Algorithm 1 (or its modification, Algorithm 2) with the computational parameters shown in Table II. For the computation of $A$'s Adversarial controller, the noise intensity was set to $\sigma = 0.1$ in (2).

Fig. 3 shows $A$'s Baseline Value function and control policy in the plane $\xi_A = \pi$. In this plane, $T$ is pointed directly away from $A$ When $T$ is relatively far from $A$ (e.g., $r > 5$), the policy resembles something akin to pure pursuit: $u_A = -\bar{u}_A \, \text{sign}(\xi_T)$, i.e., turn in the direction towards $T$. However, when $T$ is relatively close and $|\xi_T|$ is relatively small, the behavior is flipped. This corresponds to the need for $A$ to account for $T$'s WEZ at closer ranges – $A$ does not want to point directly at $T$ at close range as doing so directly aids the closure rate of $T$'s weapon (thereby putting $A$ inside $T$'s WEZ). The particular patterns of the Value and policy functions vary with $\xi_A$, but similar overall behavior is observed. The large region of zero-turn rate control in Fig. 3 corresponds to the terminal set $\mathcal{W}_A \cup \mathcal{W}_T$.

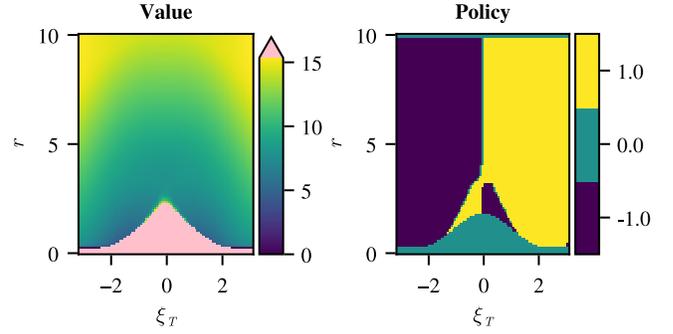

Fig. 3. Slice of the baseline Value function and policy for $\xi_A = \pi$.

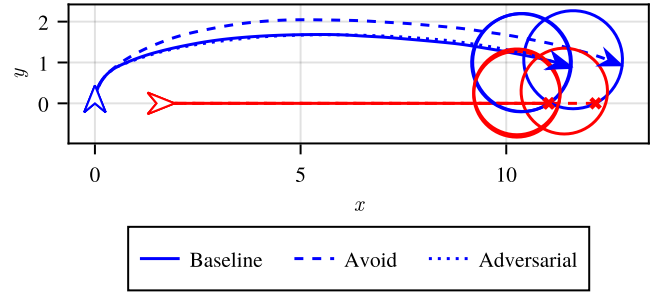

Fig. 4. Comparison of the controllers simulated against a Target moving on a straight line.

Fig. 4 shows a simulation comparing the Baseline, Avoid, and Adversarial controllers for $A$ against $T$ moving in a straight line. The open chevrons indicate initial positions while the terminal positions are denoted with filled chevrons. Each vehicle's WEZ is plotted for the terminal geometry. In all cases, the simulation ends with $\mathbf{x} \in \mathcal{W}_A \setminus \mathcal{W}_T$, that is, $T$ is in $A$'s WEZ, but $A$ is not in $T$'s WEZ. The overall behavior is also the same: $A$ turns towards $T$ up to a point and then proceeds to overtake $T$ while offset to the side, and finally turning in ahead of $T$. $A$'s Baseline and Adversarial controllers result in nearly identical trajectories while the Avoid controller maintains a larger offset as $A$ overtakes $T$ so as to avoid states for which $T$'s time-to-go would be shorter.

Indeed, for many initial conditions and $T$ controllers, $A$'s Baseline and Adversarial controller perform and behave similarly. Fig. 5 shows the results of a simulation which has the largest difference observed. For this simulation, $T$ is utilizing its Baseline controller. The lighter colors correspond to the case where $A$ uses its Baseline controller. Here, $T$ begins in a disadvantageous configuration and its optimal control is to turn hard immediately, regardless of $A$'s control. Meanwhile, it is observed that $A$'s maneuver under the Adversarial controller is slightly more aggressive than the Baseline, resulting in a slightly reduced capture time (0.91 versus 1.12 seconds).

Following is a slightly more systematic approach to comparing the Baseline and Adversarial controllers. Each pane in Fig. 6 depicts the outcomes of a simulations in which $A$ is initialized at the origin heading north, and $T$'s initial position is swept across various $x$ and $y$ coordinates along a uniform grid with the heading indicated by the chevron in the top left of the pane. Here, $A$ and $T$ both employ their Baseline controllers. Most of the space is blank due to the engagement ending in a stalemate. In those cases, the vehicles enter a configuration in which $A$ orbits $T$ who orbits a fixed point; this is partly due to the turn radius advantage of $T$. Note that $A$ has a distinct advantage in head-on configurations. Additionally, a vehicle generally wins when its opponent is placed on either side and heading towards it. Although not shown here, the number of 'losses' for $A$ does not decrease when it employs its Avoid controller.

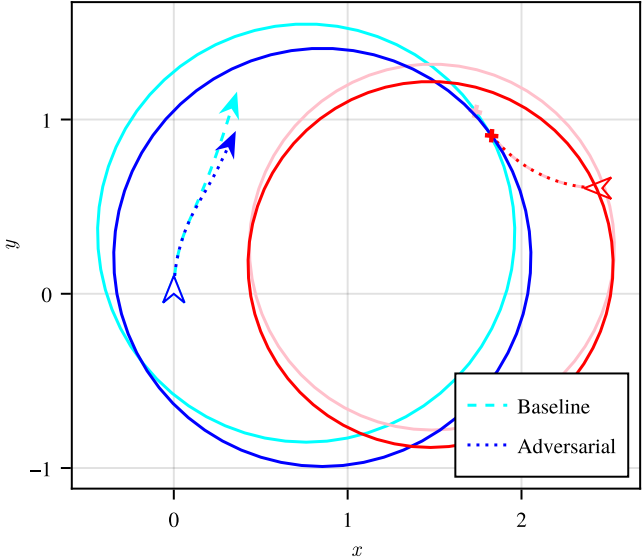

Fig. 5. Comparison of the controllers simulated against a Target using its optimal control $u_T^*(\mathbf{x})$

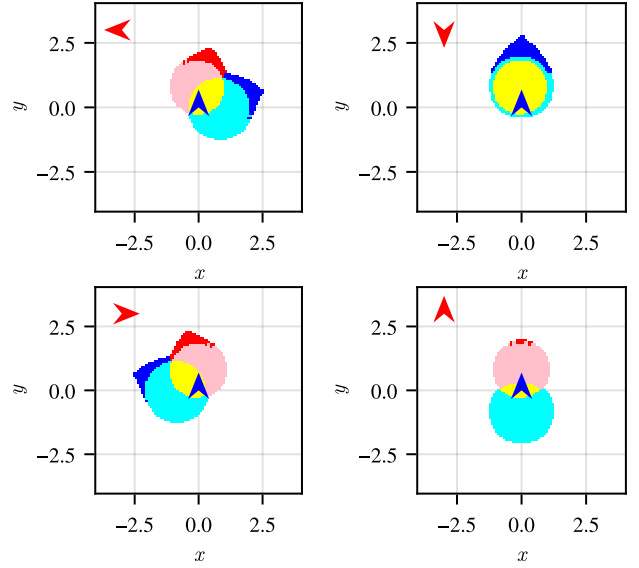

Fig. 6. Sweep of initial Target state showing the outcome of each simulation. Blank/white regions indicate initial conditions for which the simulation ends in a stalemate. Yellow regions indicate initial states in which $\mathbf{x} \in \mathcal{W}_T$ and $\mathbf{x} \in \mathcal{W}_A$. Pink and cyan regions correspond to initial states in which $\mathbf{x} \in \mathcal{W}_T$ only and $\mathbf{x} \in \mathcal{W}_A$ only, respectively. Red and blue regions correspond to initial states in which the engagement terminates in $\mathcal{W}_T$ and $\mathcal{W}_A$ respectively.

The outcomes shown in Fig. 6 are nearly identical in the case that $A$ employs its Adversarial controller. The only benefit for $A$ of employing its Adversarial controller is a very slight improvement in the capture time. Fig. 7 shows the improvement in capture times for $A$'s winning configurations for the case where $\theta_T = \pi$.

## VII. CONCLUSION

The problem of navigating a vehicle to place its weapon engagement zone onto a moving target equipped with its own weapon has been investigated. Both vehicles were modeled using Dubins dynamics, i.e., constant speed and maximum turn rate. The model corresponds to co-altitude engagements between fixed-wing aircraft. Stochastic optimal control methods were used to solve the problem; the goal of the vehicle was to get the target in its WEZ as quickly as possible. Three scenarios were solved and simulated corresponding to various assumptions regarding the target's control. In the first, the target's control was assumed to be fully unknown (i.e., just noise). The same was true of the second scenario however the vehicle sought to avoid states in which the target *could* reach it first. Finally, in the third scenario, the vehicle assumed the target uses its own optimal control. Interestingly, directly incorporating the 'knowledge' of the target's control did not appreciably improve the vehicle's performance beyond a slight improvement in capture time. Therefore, the results support the usage of the first paradigm (assuming no knowledge of the target's control or goal) for a variety of circumstances including the case where the target is adversarial. Nonetheless, the controllers designed in this paper are still based on partic-ular assumptions. Future work will focus on a full stochastic differential game treatment with the goal of obtaining robust guarantees on performance. Additionally, alternative WEZ models may be considered, along with full 3-dimensional spatial motion. Lastly, the above conclusions are made with respect to the particular parameters used in the examples; the parameter space ought to be explored, especially the noise intensity and its effects on the solutions.

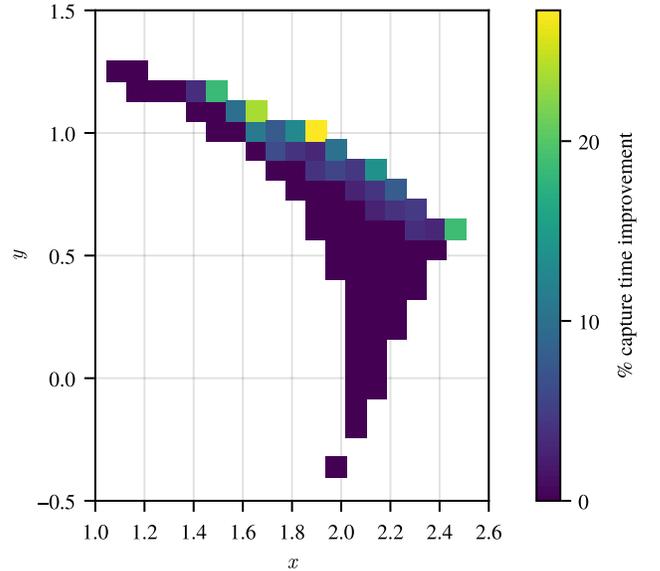

Fig. 7. Improvement in capture time for $A$'s winning configurations (with $\theta_T = \pi$) by employing the Adversarial controller (compared to Baseline).

## VIII. APPENDIX

Derivation of the $r$ increment:

$$\begin{aligned}
\mathrm{d}r &= \frac{\partial r}{\partial x_T}\mathrm{d}x_T + \frac{\partial r}{\partial y_T}\mathrm{d}y_T + \frac{\partial r}{\partial x_A}\mathrm{d}x_A + \frac{\partial r}{\partial y_A}\mathrm{d}y_A \\
&= \left(\frac{x_T - x_A}{r}v_T\cos\theta_T + \frac{y_T - y_A}{r}v_T\sin\theta_T \right.\\
&\quad \left. -\frac{x_T - x_A}{r}v_A\cos\theta_A - \frac{y_T - y_A}{r}v_A\sin\theta_A\right)\mathrm{d}t \quad (21)\\
&= (v_T\cos\lambda\cos\theta_T + v_T\sin\lambda\sin\theta_T - \\
&\quad v_A\cos\lambda\cos\theta_A - v_A\sin\lambda\sin\theta_A)\mathrm{d}t \\
&= (v_T\cos(\lambda - \theta_T) - v_A\cos(\lambda - \theta_A))\mathrm{d}t
\end{aligned}$$

which leads to (6).

Derivation of the $\xi$ increments:

$$\mathrm{d}\xi_A = \mathrm{d}\lambda - \mathrm{d}\theta_T, \quad \mathrm{d}\xi_T = \mathrm{d}\lambda - \mathrm{d}\theta_A \quad (22)$$

So we must derive $\mathrm{d}\lambda$ from

$$\begin{aligned}
\mathrm{d}\lambda &= \frac{1}{\left(\frac{y_T-y_A}{x_T-x_A}\right)^2+1}\left[-\frac{y_T-y_A}{(x_T-x_A)^2}\mathrm{d}x_T\right.\\
&\quad +\frac{1}{(x_T-x_A)^2}\mathrm{d}y_T\,\frac{y_T-y_A}{(x_T-x_A)^2}\mathrm{d}x_A \\
&\quad \left.-\frac{1}{(x_T-x_A)^2}\mathrm{d}y_A\right]\\
&= \frac{(x_T-x_A)^2}{r^2}\left[-\frac{y_T-y_A}{(x_T-x_A)^2}v_T\cos\theta_T\right.\\
&\quad +\frac{1}{(x_T-x_A)^2}v_T\sin\theta_T\,\frac{y_T-y_A}{(x_T-x_A)^2}v_A\cos\theta_A \quad (23)\\
&\quad \left.-\frac{1}{(x_T-x_A)^2}v_A\sin\theta_A\right]\mathrm{d}t\\
&= \frac{1}{r}[-v_T\sin\lambda\cos\theta_T + v_T\cos\lambda\sin\theta_T\\
&\quad + v_A\sin\lambda\cos\theta_A - v_A\cos\lambda\sin\theta_T]\mathrm{d}t\\
&= \frac{1}{r}[v_T\sin(\theta_T-\lambda) + v_A\sin(\lambda-\theta_A)]\mathrm{d}t
\end{aligned}$$

which leads to (7).